\documentclass{notices}
\usepackage{amsfonts,amssymb,amsmath,amscd}

\usepackage{graphicx} 
\usepackage[dvipsnames]{xcolor}

\usepackage{tcolorbox}
\definecolor{mycolor}{rgb}{0.122, 0.435, 0.698}

\usepackage{xurl}

\usepackage[all]{nowidow}

\title{
Shaping the Future of Mathematics in the Age of AI
}

\author{
  Johan Commelin%
  \affil{
    Johan Commelin is the director of the Mathlib Initiative and an assistant professor in the Department of Mathematics at Utrecht University. His email address is j.m.commelin@uu.nl.
    }
  \and
  Mateja Jamnik%
  \affil{
    Mateja Jamnik is a professor of artificial intelligence at University of Cambridge. Her email address is mateja.jamnik@cl.cam.ac.uk.
   }
   \and
 Rodrigo Ochigame%
 \affil{
Rodrigo Ochigame is an assistant professor of anthropology at Leiden University. Their email address is rodrigo@ochigame.org.
 }
  \and
 Lenny Taelman%
 \affil{
Lenny Taelman is a professor of mathematics at the University of Amsterdam. His email address is l.d.j.taelman@uva.nl.
 }
  \and
 Akshay Venkatesh%
 \affil{
 Akshay Venkatesh is a professor at the School of Mathematics of the Institute for Advanced Study. His email address is akshay@ias.edu.
 }
}

\date{\tiny}
\begin{document}

\maketitle

\pagebreak


Artificial Intelligence (AI) is transforming mathematics at a speed and scale that demand we reconsider the very intellectual basis of our discipline. Built on twentieth-century explorations of computability and the foundations of mathematics, modern AI has grown beyond calculation and exhaustive checking to sophisticated mathematical deduction.  The rapid emergence of these abilities creates a host of issues that our community  urgently needs to address.\looseness-1

This transformation appears across diverse and competing AI methods. A growing community of mathematicians is using the Lean proof assistant\footnote{See \url{https://leanprover-community.github.io} for the community portal about Lean and Mathlib; for associated articles, see \cite{deMoura2015,mathlibCommunity2020}.}---a formal verification system based on traditional symbolic methods---to encode and verify complex proofs in collaborative formalization projects \cite{commelin_topaz_2024}. Neuro-symbolic systems such as \mbox{AlphaProof} \cite{hubert2025alphaproof} have recently performed at medal level at the International Mathematics Olympiad, while large language models (purely neural AI) are beginning to serve as mathematical research assistants \cite{salim2025accelerating}. Though these approaches differ in method and philosophy, together they have the potential to change how we discover, verify, and organize mathematical knowledge.\looseness-1      

We---the authors of this article---organized a workshop on ``Mechanization and Mathematical Research'' in September 2025 at the Lorentz Center in Leiden, the Netherlands.%
\footnote{See \url{https://www.lorentzcenter.nl/mechanization-and-mathematical-research.html} for programme and slides, and \url{https://www.youtube.com/@mechanicalmath} for recordings of the lectures of the concluding public symposium.}
The workshop brought together mathematicians, computer scientists, philosophers and historians to look beyond the immediate technical horizon. 
The resulting discussions were electric---all  the attendees understood the stakes. 
We write primarily from our own perspective on the most important takeaways, aiming to spark a necessary conversation within the mathematical community---one that neither idealizes established practices nor assumes that technological change is inherently desirable.\looseness-1

To clarify terms for readers: by \emph{formal proofs} we mean computer-verifiable proofs in systems such as Lean, Rocq, or Isabelle. References to \emph{libraries} of formal proofs indicate shared repositories of verified statements and proofs, such as Mathlib and the Archive of Formal Proofs\footnote{\url{https://www.isa-afp.org}}. By \emph{AI} we refer to both symbolic and neural methods, including large language models. These tools differ in method but can each contribute to automating mathematical reasoning and discovery. 
Finally, we interpret \emph{mathematical community} broadly, encompassing researchers in academia and industry, as well as students and independent contributors. 

Below, we highlight five key themes that emerged from the workshop: our values in mathematics, our practice, our teaching, the technologies we engage with, and the ethical considerations arising from new intelligent systems. We will emphasize several central questions and recommendations. The future of mathematics in the age of AI need not be something that happens to us---it is something that the mathematical community should shape together.\looseness-1


\subsection*{Values in mathematics}

We were struck by the breadth of viewpoints expressed at our workshop. Some participants were enthused by the idea of computers freeing them from tedious tasks of proof writing, while others felt that these very tasks might hold some essential part of the craft of mathematics.  Some advocated for a greater role of formal proofs in research and teaching, while others defended the exact opposite.  These disagreements often highlighted an ongoing tension between professed values and the reality of mathematical practice.\looseness-1

This spectrum reflects a genuine diversity of views about basic values, tied to
mathematics' distinctive heritage as both a humanistic and a scientific discipline. 
Simplistic generalizations about the values of mathematicians
should therefore be viewed with caution.
 
Yet the current moment calls for a deliberate examination of our values in all their complexity.  \emph{Means reshape ends}: the capabilities of new technology risk dictating what counts as mathematically significant, rather than   serving our priorities. The ability to automate pieces of mathematical arguments  will change the kinds of problems that are pursued and the forms of proof that are valued;   the strengthening of rigor by means of formalization and automated proof-checking will change the expectations and standards of publication.
And  while some parts of the community may benefit from the increased commercial interest in mathematical AI, we must remain conscious of how such interest can reorient the field's internal compass. Otherwise, our discipline risks losing its intellectual autonomy.

These issues are not new; what makes them pressing is the scale and speed of recent developments.  There are no simple answers, but the workshop revealed broad agreement on a fundamental principle:\looseness-1

\begin{tcolorbox}[colback=mycolor!5!white,colframe=mycolor!75!black]
Ensure that the development and adoption of new technologies in mathematics remain firmly rooted in our own epistemic and aesthetic values---diverse as they are---rather than driven by the internal logic of those technologies or the commercial interests of their developers. 
 \looseness-1
\end{tcolorbox}

Professional societies, mathematics departments,  and individual researchers  all have important roles to play in facilitating conversations around these issues and shaping a deliberate path forward.


\subsection*{Future of mathematical practice}

We do not know if, when, or to what extent AI systems will surpass humans in their ability to perform clearly defined mathematical tasks such as generating formal proofs. The discussions at the Lorentz Center, however, underscored that it is well past time to grapple with the implications of these possibilities.

As a community, we must now ask: in such scenarios, what will qualify as mathematical discovery? Will it be the generation of a formal proof, the formulation of an important conjecture, the construction of an explanatory argument, or something else? What will the role of human intuition be?  How can we nurture the next generation of mathematicians when the very nature of our practice is in flux? Such questions help us frame our responsibility:

 \begin{tcolorbox}[colback=mycolor!5!white,colframe=mycolor!75!black]
Engage collectively---and alongside our students---with the question of how our relationship to mathematics  will change if computers can perform many of the day-to-day tasks of current mathematical research.\looseness-1
\end{tcolorbox}


\subsection*{Teaching and education}

As with the introduction of handheld calculators and other past technologies, AI may change what we regard as fundamental mathematical competencies and the very purpose of mathematical education. If computers can perform precisely defined mathematical tasks at scale, then mathematical fluency may no longer lie in executing those tasks oneself. Instead, it may lie in understanding which tasks to perform and why, in connecting them to larger questions, and in cultivating the judgment and insight that give mathematics its meaning. This shift also raises questions about assessment and examinations: how can we evaluate students’ understanding and skills in a landscape where AI can perform routine calculations and proofs? To make the mathematical community more resilient to an unpredictable future, we propose that the community:

\begin{tcolorbox}[colback=mycolor!5!white,colframe=mycolor!75!black]
    Review undergraduate and graduate mathematical curricula to emphasize a broader range of skills, including posing problems, communicating ideas, and critiquing purported logical arguments. 
   \looseness-1
\end{tcolorbox}

Students should have a say in this deliberation: their future is at stake.


\subsection*{Technology}

How can mathematicians contribute directly to the creation of AI tools for mathematics, so that we have more of a voice in their development?   

One way this is already happening is through formal proof libraries.  These libraries have succeeded in involving a large community of professional and non-professional mathematicians. To ensure long-term health of these libraries and their alignment with academic values, we must now address critical questions of governance, learning from both the successes and the failures of open-source projects. In particular, we must establish clear licensing frameworks and contributor norms that protect these libraries from unrestricted commercial use without attribution or respect for the rights of authors.
 
Another proposal being pursued by multiple groups is the creation of community-owned benchmarks. Good benchmarks can signal community values to those developing AI systems, while also helping mathematicians understand realistic use cases for, and limitations of, the technology. For example, benchmarks that focus on open-ended problems can counter overstated claims based on narrow task suites.\looseness-1

At present, the most advanced AI systems for mathematical reasoning are proprietary and originate in for-profit labs, yet open-source and open-weights alternatives are emerging (e.g., \cite{godelprover, kimina}). These different implementations naturally compete, but they also create opportunities for synergy, both within academia and between academia and industry: shared learning, cross-validation, and complementary approaches can strengthen the community’s understanding and development of AI tools. 

The mathematical community thus faces an ambitious goal:

\begin{tcolorbox}[colback=mycolor!5!white,colframe=mycolor!75!black]
Build academically oriented technological infrastructure, with open-source code, transparent decision-making, and community oversight, in ways that reduce dependence on any single commercial actor.\looseness-1
\end{tcolorbox}

Besides helping to maintain intellectual independence,  such infrastructure would allow for more transparency around issues of cost, computational resources, training data, and environmental impact,  and would help to align technological development with the community’s values and priorities.


\subsection*{Norms, ethics, and governance}

The integration of AI into mathematical research raises significant ethical and governance concerns. The training of such systems uses the work of the mathematical community in an opaque way that threatens established norms of attribution and credit.   AI-generated papers are starting to put still more pressure on a publication system already strained to the breaking point; and how should formal proofs be integrated into this? Mechanization may also create new forms of inequality, in which disparities in access to computational resources distort the competitive dynamics of mathematical research. Moreover, the substantial energy demands of AI systems raise important environmental concerns. 

We thus propose that the community:

\begin{tcolorbox}[colback=mycolor!5!white,colframe=mycolor!75!black]
Develop and maintain a living statement of ethical principles to guide academic mathematicians and mathematical institutions in their interactions with AI systems and developers.\looseness-1
\end{tcolorbox} 

 Such a statement is not intended to be binding or exhaustive, but rather a collective articulation of shared values and ethical practices, established through wide consultation
 (compare with \cite{asilomar1975}). 
It may address norms about attribution and citation, disclosure of computational resources and training data, guidelines on the licensing of preprints and formal proof libraries, and a code of conduct for AI use. Maintaining such a statement requires ongoing governance, ideally under the leadership of professional societies, to adapt to new challenges. At the workshop, a group of participants started developing such a
statement.


\subsection*{Please engage!}

The torrent of AI-related news may provoke skepticism, or, conversely, a sense of resignation. Yet we exhort everyone reading this to engage with the topic of AI and mathematics.
Learn about its technical aspects and capabilities. 
Discuss it with your colleagues, your students, and your professional organizations. Reflect on what you value in your practice of mathematics, and take collective action.

 We are at a crossroads. As we see it, the nature of the discipline is being renegotiated. 
 By engaging with these questions, you are helping to ensure that this future will reflect our collective values and aspirations.\looseness-1

\bibliography{leiden-ref}
\end{document}